\documentclass{gtart_h}

\def\ifplaintex{\expandafter\ifx\csname documentclass\endcsname\relax}

\def\gtp{{\mathsurround=0pt\it $\cal G\mskip-2mu$eometry \&\ 
$\cal T\!\!$opology $\cal P\!$ublications}}  

\def\recd{{\small Received:\qua\receiveddate\ifx\reviseddate\relax
\else\qquad Revised:\qua\reviseddate\fi\par}} 

\def\lognumber#1{\def\thelognumber{#1}}
\def\volumenumber#1{\def\thevolumenumber{#1}}
\def\volumeyear#1{\def\thevolumeyear{#1}}
\def\papernumber#1{\def\thepapernumber{#1}}
\def\pagenumbers#1#2{\def\startpage{#1}\def\finishpage{#2}}
\def\published#1{\def\publishdate{#1}}

\def\received#1{\def\receiveddate{#1}}
\def\revised#1{\def\reviseddate{#1}}
\def\accepted#1{\def\accepteddate{#1}}

\long\def\asciiabstract#1{\long\def\theasciiabstract{#1}}

\let\\\par\let\thelognumber\relax\let\thevolumenumber\relax
\let\thepapernumber\relax\let\thevolumeyear\relax\let\startpage\relax
\let\finishpage\relax\let\publishdate\relax\let\receiveddate\relax
\let\reviseddate\relax\let\accepteddate\relax\let\theasciititle\relax
\let\theasciiauthors\relax
\let\theasciiabstract\relax

\let\theasciiemail\relax

\ifplaintex
\font\logobig=cmssbx10 scaled 3836
\font\logomed=cmssbx10 scaled 2557
\else
\font\logobig=cmssbx10 scaled 4200
\font\logomed=cmssbx10 scaled 2800
\fi

\long\def\makeagttitle{   
\count0=\startpage
\agt\hfill      
\hbox to 45truept{\vbox to 0pt{\vglue -13truept{\logomed A\kern -.37em{\logobig 
T}\kern -.38em G}\vss}\hss}
\break
{\small Volume \thevolumenumber\ (\thevolumeyear)
\startpage--\finishpage\nl
Published: \publishdate}

\vglue .25truein

{\parskip=0pt\leftskip 0pt plus
1fil\def\\{\par\smallskip}{\Large\bf\thetitle}\par\medskip} \vglue
0.05truein

{\parskip=0pt\leftskip 0pt plus 1fil\def\\{\par}{\sc\theauthors}
\par\medskip}%
 
\vglue 0.03truein 

{\small\leftskip 25truept\rightskip 25truept{\bf Abstract}\stdspace\theabstract

{\bf AMS Classification}\stdspace\theprimaryclass
\ifx\thesecondaryclass\relax\else; \thesecondaryclass\fi\par
{\bf Keywords}\stdspace \thekeywords\par}\vglue 7truept

}   

\ifplaintex
\font\phead=cmsl9 scaled 950
\font\pnum=cmbx10 scaled 913
\font\pfoot=cmsl9 scaled 950
\headline{\vbox to 0pt{\vskip -4.5mm\line{\small\phead\ifnum
\count0=\startpage ISSN 1472-2739 (on-line) 1472-2747 (printed)
\hfill {\pnum\folio}\else\ifodd\count0\def\\{ }%
\ifx\theshorttitle\relax\thetitle\else\theshorttitle\fi\hfill{\pnum\folio}
\else\def\\{ and }{\pnum\folio}\hfill\ifx\theshortauthors\relax\theauthors
\else\theshortauthors\fi\fi\fi}\vss}}
\footline{\vbox to 0pt{\vglue 0mm\line{\small\pfoot\ifnum\count0=\startpage
\copyright\ \gtp\hfill\else
\agt, Volume \thevolumenumber\ (\thevolumeyear)\hfill\fi}\vss}}
\else
\font=cmsl9 scaled 1050
\font\lnum=cmbx10 
\font=cmsl9 scaled 1050
\makeatletter
\def\@oddhead{{\small\lhead\ifnum\count0=\startpage ISSN 1472-2739 
(on-line) 1472-2747 (printed)\hfill {\lnum\number\count0}\else\ifodd\count0
\def\\{ }\ifx\theshorttitle\relax \thetitle \else\theshorttitle\fi\hfill
{\lnum\number\count0}\else\def\\{ and }{\lnum\number\count0}
\hfill\ifx\theshortauthors\relax 
\theauthors\else\theshortauthors\fi\fi\fi}}\def\@evenhead{\@oddhead}
\def\@oddfoot{\small\lfoot\ifnum\count0=\startpage\copyright\ \gtp\hfill\else
\agt, Volume \thevolumenumber\ (\thevolumeyear)\hfill\fi}
\def\@evenfoot{\@oddfoot}
\makeatother
\fi
\let\maketitlepage\makeagttitle
\let\makeshorttitle\maketitlepage
\let\maketitle\maketitlepage

\newwrite\gtoutfile
\long\gdef\makeheadfile{  
{\def\\{, }\def\s{ }
\immediate\openout\gtoutfile head.xxx
\immediate\write\gtoutfile{Proxy-for: \ifx\theasciiauthors\relax
\theauthors\else\theasciiauthors\fi\s<\ifx\theasciiemail\relax\theemail\else\theasciiemail\fi>}
\immediate\write\gtoutfile{\noexpand\\}
\immediate\write\gtoutfile{Authors: \ifx\theasciiauthors\relax
\theauthors\else\theasciiauthors\fi}
{\def\\{ }\immediate\write\gtoutfile{Title: \ifx\theasciititle\relax
\thetitle\else\theasciititle\fi}}
\immediate\write\gtoutfile{Subj-class: GT or SG, GR etc}
\immediate\write\gtoutfile{MSC-class: \theprimaryclass\ifx\thesecondaryclass\relax\else, \thesecondaryclass\fi}
\immediate\write\gtoutfile{Journal-ref: Algebr. Geom. Topol. \thevolumenumber\s
(\thevolumeyear) \startpage-\finishpage}
\immediate\write\gtoutfile{Comments: Published by Algebraic and
Geometric Topology at}
\immediate\write\gtoutfile{\s\s\s  http://www.maths.warwick.ac.uk/agt/AGTVol\thevolumenumber/agt-\thevolumenumber-\thepapernumber.abs.html}
\immediate\write\gtoutfile{\noexpand\\}
\immediate\write\gtoutfile{}
\ifx\theasciiabstract\relax
\immediate\write\gtoutfile{\theabstract}\else
\immediate\write\gtoutfile{\theasciiabstract}\fi
\immediate\write\gtoutfile{}
\immediate\write\gtoutfile{\noexpand\\}
\immediate\write\gtoutfile{}
\immediate\closeout\gtoutfile}}  

\def\maketitlepage{\makeagttitle\makeheadfile}
\let\makeshorttitle\maketitlepage
\let\maketitle\maketitlepage

\lognumber{27}
\volumenumber{4}
\volumeyear{2004}
\papernumber{27}
\published{5 November 2003}
\pagenumbers{595}{602}
\received{2 January 2003}
\revised{25 March 2004}
\accepted{12 July 2004}

\usepackage{amssymb, amsmath}

\theoremstyle{plain}
\newtheorem{theorem}{Theorem}[section]

\newtheorem{lemma}[theorem]{Lemma}

\newtheorem{proposition}[theorem]{Proposition}

\theoremstyle{definition}

\newtheorem*{acknowledgements}{Acknowledgements}

\newcommand{\Z}{{\mathbb{Z}}}
\newcommand{\F}{{\mathbb{F}_2}}
\newcommand{\R}{{\mathbb{R}}}

\begin{document}
\title{Commutators and squares in free groups}
\author{Sucharit Sarkar}
\email{bmat0212@isibang.ac.in}
\address{Stat-Math Unit,
Indian Statistical Institute\\Bangalore, India}

\begin{abstract}
Let ${\mathbb{F}_{2}}$ be the free group generated by $x$ and $y$. In this
article, we prove that the commutator of $x^{m}$ and $y^{n}$ is a product of
two squares if and only if $mn$ is even. We also show using 
topological methods that there are
infinitely many obstructions for an element in ${\mathbb{F}_{2}}$ to be a
product of two squares.
\end{abstract}
\asciiabstract{%
Let F_2 be the free group generated by x and y. In this article, we
prove that the commutator of x^m and y^n is a product of two squares
if and only if mn is even. We also show using topological methods that
there are infinitely many obstructions for an element in F_2 to be a
product of two squares.}


\maketitle

\section{Introduction}

In any group $G$, the commutator of two elements $g$ and $h$ is a
product of three squares, namely,
$$ghg^{-1}h^{-1}=g^{2}\left( g^{-1}h\right) ^{2}h^{-2}.$$ Let $\left[
g,h\right] $ denote the commutator of $g$ and $h$. It is natural to
ask whether $\left[ g,h\right] $ can be written as a product two
squares. Since the subgroup generated by $g$ and $h$ in $G$ is a
quotient of the free group on two generators, the answer would be in
the affirmative if we knew that the commutator of the generators of
the free group on two generators can be written as a product of two
squares. However, if ${\mathbb{F}_{2}}$ is the free group on two
generators $x$ and $y$, a theorem of Lyndon and Newman~\cite{LN} states
that the commutator $[x,y]$ is not a product of two squares. Here, we
give the following generalisation of their theorem. Further, the
method of our proof extends to give infinitely many obstructions to an
element being the product of two squares.

\begin{theorem}\label{Lyn} 
$[x^{m},y^{n}]$ is a product of two squares in
${\mathbb{F}_{2}}=\langle x,y\rangle$ if and only if $mn$ is even.
\end{theorem}

In the case when $m=n=1$ this gives the theorem of Lyndon and
Newman. Our methods also give a proof of the following theorem of
Akhavan-Malayeri~\cite{Ak}.

\begin{theorem}[Akhavan-Malayeri]\label{AM}
$[x,y]^{2n+1}$ is not a product of two squares in
${\mathbb{F}_{2}}=\langle x,y\rangle$.
\end{theorem}

We first reformulate the question in terms of products of conjugate elements.

\begin{lemma}\label{conjsq}
An element $g\in {\mathbb{F}_{2}}$ can be expressed as a product $%
g=a^{2}b^{2}$ of two squares if and only if $g$ can be expressed as a
product $g=c(d^{-1}cd)$ of two conjugate elements.
\end{lemma}

\begin{proof}
This follows immediately from the relation $a^2b^2=(ab)(b^{-1}(ab)b)$.
\end{proof}

The following statement underlines the importance of the previous lemma. 

\begin{lemma}
Suppose $g\in [{\mathbb{F}_2},{\mathbb{F}_2}]$ is of the form $g=ab$
with $a$ and $b$ conjugate in $\F$. Then $a,b\in
[{\mathbb{F}_2},{\mathbb{F}_2}]$.
\end{lemma}

\begin{proof}
Let $g\mapsto\bar g$ be the abelianisation homomorphism
$\F\to\Z^2$. Then as $g\in [\F,\F]$, $0=\bar g=\bar a + \bar b=2\bar a$
where the last equality holds as $a$ and $b$ are conjugate. Thus,
$\bar a=\bar b=0$ which implies that $a,b\in [\F,\F]$.
\end{proof}

The heart of our proof lies in constructing a group homomorphism $\varphi \co [{%
\mathbb{F}_{2}},{\mathbb{F}_{2}}]\rightarrow {\mathbb{Z}}$ which is
invariant under conjugacy action of ${\mathbb{F}_{2}}$, i.e., $\varphi
(a)=\varphi (g^{-1}ag)$ for all $a\in \lbrack {\mathbb{F}_{2}},{\mathbb{F}%
_{2}}]$ and $g\in {\mathbb{F}_{2}}$. The construction of our homomorphism is
topological.

If an element $g\in [\F,\F]$ is a product of two squares, then we
shall see that $\varphi(g)$ is even. This gives a criterion to decide
whether $g$ is the product of two squares. The theorem of Lyndon and
Newman follows from this. We shall extend this to stronger criteria
for an element in $[\F,\F]$ to be a product of two squares.

\section{The homomorphism $\protect\varphi$}

We construct a homomorphism
$\varphi\co [{\mathbb{F}_{2}},{\mathbb{F}_{2}}]\rightarrow
{\mathbb{Z}}$ which is invariant under the conjugacy action of
${\mathbb{F}_{2}}$. The theorem of Lyndon and Newman follows from the
properties of this homomorphism.

Let $K$ be the wedge of two circles. Then $\pi_1(K)=\F$. Let $\tilde
K$ be the universal abelian cover of $K$. This is the cover
corresponding to the subgroup $[\F,\F]$ of $\F$.  We can identify
$\tilde K$ with a subcomplex of $\R^2$ whose vertices are $\Z^2$ and
edges join $(i,j)$ either to $(i+1,j)$ or to $(i,j+1)$.

Let $X$ and $Y$ denote the edges from $(0,0)$ to $(1,0)$ and $(0,1)$
respectively. Under the group $\Z^2$ of deck transformations each edge
is the image of $X$ or $Y$. Denoting the group of deck transformations
multiplicatively and taking $x$ and $y$ to be the generators, we see
that all edges are of the form $x^iy^jX$ or $x^iy^jY$ with $i,j\in\Z$.

Consider the simplicial homology of $\tilde K$. We shall first define
a homomorphism $\theta\co Z_1(\tilde K)\to\Z$ from the $1$-cycles to the
integers which is invariant under deck transformations. The
homomorphism $\varphi$ will be defined in terms of this.

Observe that, using the above notation, the simplicial chains of
$\tilde K$ are of the form $\alpha=P_\alpha(x,y)X + Q_\alpha(x,y) Y$
with $P_{\alpha}(x,y)$ and $Q_{\alpha}(x,y)$ Laurent polynomials. Let
$f_\alpha(y)=P_\alpha(1,y)$. The action by a deck transformation $x^ky^l$
takes $\alpha=P_\alpha(x,y)X + Q_\alpha(x,y) Y$ to
$x^ky^l\alpha=x^ky^lP_\alpha(x,y)X + x^ky^lQ_\alpha(x,y) Y$.

Further, for a cycle $\alpha$, it is easy to see that
$f_\alpha(1)=P_\alpha(1,1)=0$. We define the homomorphism $\theta$ by
$$\theta(\alpha)=f'_\alpha(1).$$ 

\begin{lemma}
$\theta\co Z_1(\tilde K)=H_1(\tilde K)\to\Z$ is invariant under the
group of deck transformations.
\end{lemma}
\begin{proof}
It suffices to show that
$\theta(x\alpha)=\theta(\alpha)=\theta(y\alpha)$. The first equality is
obvious as $f_{x\alpha}=f_\alpha$. The second follows as
$$f_{y\alpha}'(1)=(yf_\alpha)'(1)=f_\alpha(1)
+f'_\alpha(1)=f'_\alpha(1).$$ Here we used the fact that for a cycle
$f_\alpha(1)=0$.
\end{proof}

As $\tilde K$ is a 1-complex, this can be viewed as a homomorphism,
also denoted $\theta$, from $H_1(\tilde K)$ to $\Z$. Using this, we
define $\varphi\co [\F,\F]\to\Z$. Namely, given a curve $\gamma$ in $K$
representing an element $g\in [\F,\F]$, take a lift $\tilde\gamma$ of
$\gamma$ to $\tilde K$. This represents an element in homology, and we
let $\varphi(g)=\theta(\tilde\gamma)$. This is independent of the lift
chosen as different lifts are related by deck transformation, and
$\theta$ is invariant under deck transformations. 

We need some properties of $\varphi$.

\begin{lemma}
$\varphi(g)=\varphi(hgh^{-1})$ for all $g\in [\F,\F]$ and $h\in \F$.
\end{lemma}
\begin{proof}
The elements of $g$ and $hgh^{-1}$ of $[\F,\F]$ can be represented by
curves that have lifts in $\tilde K$ that differ by deck
transformations. As deck transformations leave $\theta$
invariant, $\varphi(g)=\varphi(hgh^{-1})$.
\end{proof}
 
\begin{lemma}
Suppose $g\in [\F,\F]$ is a product of two squares. Then $\varphi(g)$
is even.
\end{lemma}
\begin{proof}
By Lemma~\ref{conjsq} we can write $g=hk$ with $h$ and $k$ elements of
$[\F,\F]$ that are conjugate in $\F$. We have
$\varphi(g)=\varphi(h)+\varphi(k)=2\varphi(h)$ as $\varphi$ is
invariant under conjugation.
\end{proof}

Consider the lift to $\tilde K$ of a curve $\gamma$ representing
$[x,y]$ in $K$. This lift gives the cycle $\alpha=(1-y)X+(x-1)Y$ and,
hence $\varphi([x,y])=f'_\alpha(1)=-1$. Thus, $\varphi$ is
non-trivial. Moreover, by the above lemma $[x,y]$ is not a product of
two squares.

We now prove our extension of the result of Lyndon and Newman, one
half of which is an extension of the above argument.

\begin{proof}[Proof of Theorem~\ref{Lyn}]
Take a curve $\gamma$ representing $[x^m,y^n]$. We see that the lift
of this starting at the origin gives the chain 
$$\alpha=(1+x+\dots x^{m-1})(1-y^n)X+(1+y+\dots+y^{n-1})(x^m-1)Y$$ and
hence $\varphi([x^m,y^n])=-mn$. Hence if $[x^m,y^n]$ is the product of two
squares, $mn$ is even.

Conversely, if $mn$ is even, we assume without loss of generality that
$m$ is even. Then
$[x^m,y^n]=(x^{\frac{m}{2}})^2(y^nx^{\frac{-m}{2}}y^{-n})^2$.
\end{proof}

The same methods also yield a proof of Theorem~\ref{AM} of
Akhavan-Malayeri. We thank the referee for pointing this out.

\begin{proof}[Proof of Theorem~\ref{AM}]
By the above, $\varphi([x,y]^{2n+1})=-(2n+1)$ as $\varphi$ is a
homomorphism. The result follows.
\end{proof}

We have constructed a homomorphism $\varphi $ on
$[\mathbb{F}_{2},\mathbb{F}_{2}]$ that is invariant under
conjugacy. Now $[{\mathbb{F}_{2}},{\mathbb{F}_{2}}]$ is the smallest
normal subgroup of $\mathbb{F}_{2}$ containing $ \left[ x,y\right]
$. Therefore, any such homomorphism is determined by its value on
$[x,y]$. In particular, if we make the analogous construction taking
$g(x)=Q_{\alpha}(x,1)$ in place of $f(y)$ and define $\psi (\alpha
)=g^{\prime }(1)$, we have $\psi =-\varphi $ as we can see by
evaluating on $[x,y]$.

\section{Further criteria}

So far we have one criterion for an element in
$[{\mathbb{F}_2},{\mathbb{F}_2} ]$ to be a product of two squares,
namely, if $g$ is a product of two squares then $\varphi(g)$ is even.
If $\varphi(g)=0$ (which implies $\psi(g)=0$), there are other
criteria. These are obtained by constructing homomorphisms $\varphi_2$
and $\psi_2$ on appropriate subgroups of $G=[\F,\F]$ which are even on
elements that are products of two squares.

We first need some lemmas. Let $G_{1}$ denote the kernel of $\varphi
$. Since $\varphi$ is conjugacy invariant, $G_{1}$ is a normal
subgroup of ${\mathbb{F}_{2}}$. Let $H_{1}= ker (\psi) = G_{1}$.

\begin{lemma}
\label{ingp} Suppose $\varphi (g)=0$ and $g=ab$ with $a$ and $b$ conjugate.
Then $a,b\in G_{1}$.
\end{lemma}

\begin{proof}
As $\varphi$ is conjugacy invariant, $0=\varphi(g)=2\varphi(a)$. 
\end{proof}

We need an elementary property of derivatives of polynomials.

\begin{lemma}
Let $f(t)$ be a Laurent polynomial with integer coefficients. Then
$f^k(t)$ is divisible by $k!$.
\end{lemma}

\begin{proof}
The $k$th derivative of $t^n$ is divisible by $n(n-1)\cdots (n-k+1)$
which in turn is divisible by $k!$. The result follows.
\end{proof}

Now we can define two homomorphisms $\varphi_2$ and $\psi_2$ from
$H_1=G_1$ to ${\mathbb{Z}}$ by
$\varphi_2(g)=f_\alpha^{\prime\prime}(1)/2$ and $
\psi_2=g^{\prime\prime}_\alpha(1)/2$ with $\alpha$ as before. The
proof that these are well defined and conjugacy invariant is exactly
as in the previous section.

Continuing in this manner, we let $G_{2}=ker(\varphi _{2})$ and $
H_{2}=ker(\psi _{2})$. We inductively define groups $H_{k}$ and
$G_{k}$ and homomorphisms $\varphi_k\co G_{k-1}\to \Z$ and $\psi_k\co H_{k-1}\to
\Z$. Namely, let $\varphi _{k}(\alpha )=f_{\alpha }^{(k)}(1)/k!$ and
$\psi _{k}(\alpha )=g_{\alpha }^{(k)}(1)/k!$ and define
$G_{k}=ker(\varphi _{k})$ and $H_{k}=ker(\psi _{k})$.

As in the previous section, we deduce the following properties of the
homomorphisms $\varphi_k$ and $\psi_k$.

\begin{lemma}
The homomorphisms $\varphi_k$ and $\psi_k$ are invariant under the
action of $\F$ by conjugation.
\end{lemma}

\begin{lemma}
If $g\in G_k$ (respectively $g\in H_k$) is a product of two
squares, then $\varphi_k(g)$ (respectively $\psi_k$) is even.
\end{lemma}

Thus, we have infinitely many obstructions to an element being the
product of two squares. More precisely, let $g\in [\F,\F]$ be an
element. We evaluate $\varphi(g)=-\psi(g)$. There are three
possibilities: $\varphi(g)=0$, $\varphi(g)$ is even or $\varphi(g)$ is
odd. In case $\varphi(g)$ is odd, we know that $g$ is not a square. If
it is even but non-zero, we cannot deduce any further obstructions. In
the case when $\varphi(g)=0$, we have additional homomorphisms
$\varphi_2$ and $\psi_2$ which can be applied to $g$ to get an odd or
even number. This process can be continued inductively.

\section{Some examples}

We have constructed in the previous section two sequences of
obstructions to an element $g\in [\F,\F]$ being a product of two
squares, based on the homomorphisms $\varphi_k$ and $\psi_k$. We shall
show that all these are non-trivial in the sense that there are
elements for which the first $k-1$ homomorphisms vanish and the $k$th is
odd.

On the other hand, in the case when $f(y)=0$ and $g(x)=0$ as
polynomials, all our homomorphisms vanish. We shall construct examples
where this happens.

Our examples are based on the observation that the associations
$\gamma\to P_\gamma$ and $\gamma\to Q_\gamma$ are module homomorphisms
over the ring of Laurent polynomials in variables $x$ and $y$, from
$H_1(\tilde K)$ to Laurent polynomials. Further, we have a surjection
from the commutator subgroup $[\F,\F]=\pi_1(\tilde K)$ to its
abelianisation $H_1(\tilde K)$.

\begin{proposition}
For any $k>1$, there is an element $g\in [\F,\F]$ such that
$\varphi_j(g)=0$ for all $j<k$, $\varphi_k(g)=-1$ and
$\psi_j(g)=0$ for all $j$.
\end{proposition}
\begin{proof}
Let $g$ be an element whose image in $H_1(\tilde K)$ is
$\gamma=(y-1)^{k-1}\alpha$, where $\alpha$ denotes the class of
$[x,y]$ in $H_1(\tilde K)$. As $P_\alpha=-(y-1)$, and the associations
$\gamma\to P_\gamma$ and $\gamma\to Q_\gamma$ are module
homomorphisms, it follows that $P_{\gamma}=-(y-1)^k$ and
$Q_{\gamma}=(y-1)^{k-1}(x-1).$

Now $f_\gamma(y)=P_\gamma(1,y)=-(y-1)^k$ and
$g_\gamma(x)=Q_\gamma(x,1)=0$ and hence $\varphi_j(g)=0$ for all
$j<k$, $\varphi_k(g)=-1$ and $\psi_j(g)=0$ for all $j$.
\end{proof}

Recall that we consider succesively the homorphisms $\varphi_k$
(and $\psi_k$), with $\varphi_{k+1}$ defined if $\varphi_k$
vanishes. The first non-zero $\varphi_k$ gives a criterion for an
element being the product of two squares. We see in the next example
that there are elements in $[\F,\F]$ for which all the $\varphi_k$
vanish.

\begin{proposition}
There is an element $g\in [\F,\F]$ such that
$\varphi_j(g)=0 = \psi_j(g)$, for all $j$.
\end{proposition}
\begin{proof}
We take $g$ whose image in $H_1(\tilde K)$ is
$\gamma=(x-1)(y-1)\alpha$. Then, using the notation of the above proposition,
we see that $f_\gamma(y)=0$ and $g_\gamma(x)=0$, and hence
$\varphi_j(g)=0 = \psi_j(g)$, for all $j$.
\end{proof}

\section{Factorisation and another criterion}

The above example suggests a variant of our criteria.

\begin{proposition}
Let $g\in [\F,\F]$ be an element with $\alpha$ the corresponding
cycle in $H_1(\tilde K)$. Suppose $(x-1)^k(y-1)^l$ divides
$P_\alpha(x,y)$ with quotient $h(x,y)$. If $g$ is a product of
two squares then $h(1,1)$ is even.
\end{proposition}
\begin{proof}
We have seen that if $g$ is the product of two squares, it is the
product of two conjugates $c$ and $d$.  Let $Ab$ denote the
abelianisation map $[\F,\F]\to H_1(\tilde K)$. Then
$Ab(g)=Ab(c)+Ab(d)$. As $c$ and $d$ are conjugate (in $\F$), they
differ by a deck transformation of $\tilde K$. Equivalently, 
$c$ and $d$ differ by the action by conjugation of the abelianisation
$\F/[\F,\F]=\Z^2$ on the commutator subgroup $[\F,\F]$. This action has been
identified with multiplication by Laurent polynomials.

Hence, for some integers $m,n\in\Z$, we have $Ab(d)=x^my^n Ab(c)$. It follows
that $Ab(g)=(x^my^n+1)Ab(c)$ and hence $(1+x^my^n)$ divides
$P_{\alpha}$.

Now the ring of Laurent polynomials over $\Z$ in $x$ and $y$ is a
unique factorisation domain, and $x-1$ and $y-1$ are prime
elements. Further, they do not divide $E(x,y)=(1+x^my^m)$ as
$E(1,1)\neq 0$. Thus, $E(x,y)$ divides $h(x,y)$. As $E(1,1)=2$,
$h(1,1)$ is even.
\end{proof}

\begin{acknowledgements}
I would like to thank S. Gadgil, S. P. Inamdar and B. Sury for their
advice and encouragement, which also enabled me to refine what I had
first proved. I thank the referee for several helpful comments.
\end{acknowledgements}

\newpage\bibliographystyle{amsplain}

\Addresses\recd

\end{document}